\documentclass[preprint,fleqn,11pt]{elsarticle}
\usepackage{amssymb}
\biboptions{comma,square}

\journal{Journal of Computational Physics}
\usepackage{amssymb,amsmath}
\usepackage{graphics}
\usepackage{graphicx}
\usepackage{color}

\newcommand{\dl}[1]{{\bf Theorem{#1.}}}
\newcommand{\yl}[1]{{\bf Lemma{#1.}}}

\newcommand{\zb}{\textbf{\qquad$\Box$}}
\newcommand{\la}[1]{\label{#1}}
\newcommand{\rf}[1]{(\ref{#1})}

\newcommand{\zm}{{\bf Proof.}}
\newcommand{\commentout}[1]{{}}

\begin{document}
\thispagestyle{empty} \setcounter{page}{1}

\begin{frontmatter}

\title{The Weak Galerkin Finite Element Method for the Transport-Reaction Equation}

\author{Tie Zhang\footnote{Corresponding author at: Department of
Mathematics, Northeastern University, Shenyang 110004, China. {\em
E-mail address} : ztmath@163.com (T. Zhang). Tel \& Fax:
+86-024-83680949. },\quad Shangyou Zhang$^\dag$
}
\address{$^*$Department of Mathematics, Northeastern University, Shenyang 110004, China\\
$^\dag$Department of Mathematical Sciences, University of Delaware, Delaware, 19716, USA}

\begin{abstract}
We present and analyze a weak Galerkin finite element method for solving the transport-reaction equation in $d$ space dimensions. This method is highly flexible by allowing the use of discontinuous finite element on general meshes consisting of arbitrary polygon/polyhedra. We derive the $L_2$-error estimate of $O(h^{k+\frac{1}{2}})$-order for the discrete solution when the $k$th-order polynomials are used for $k\geq 0$. Moreover, for a special class of  meshes, we also obtain the optimal error estimate of $O(h^{k+1})$-order in the $L_2$-norm. A derivative recovery formula is presented to approximate the convection directional derivative and the corresponding superconvergence estimate is given. Numerical examples on compatible and non-compatible meshes are provided to show the effectiveness of this weak Galerkin method.
\end{abstract}

\begin{keyword}
Weak Galerkin method; transport-reaction equation;
 optimal error estimate; derivative approximation; superconvergence

\MSC 65M60, 65N30, 65N12
\end{keyword}
\end{frontmatter}
\section{Introduction}
\setcounter{section}{1}\setcounter{equation}{0} In this work, we study a weak Galerkin finite element (WG) method for solving the transport-reaction equation
\begin{eqnarray}
\left\{\begin{array}{ll}
\nabla\cdot(\beta u)+\alpha u=f,\;in\;\;\Omega,\\
u=g,\;\;on\;\;\partial\Omega_-,\label{1.1}
\end{array}
 \right.
\end{eqnarray}
where $\nabla\cdot$ is the divergence operator, $\Omega\subset R^d$ ($d=2,3$) is a polygonal ($d=2$) or polyhedral ($d=3$) domain with boundary $\partial\Omega$, $\partial\Omega_-=\{x\in\partial\Omega: \beta(x)\cdot{\bf n}(x)<0\}$ and ${\bf n}(x)$ is the outward unit normal vector at the point $x\in\partial\Omega$. \textcolor[rgb]{0.00,0.00,0.00}{This equation is commonly used to model the neutron transport process.} In this case, function $u(x)$ represents the flux of neutrons at point $x$ in the angular direction $\beta$, $\alpha$ is the nuclear cross section and the source term $f$ stands for the scattering and the fission. The boundary condition means that neutrons are entering the system from inflow boundary $\partial\Omega_-$.

Equation \rf{1.1} is a typical first order hyperbolic equation. It is well known that just for problem \rf{1.1}, the original discontinuous Galerkin finite element (DG) method was presented by Reed and Hill in 1973 \cite{Reed}. Since then, many research works have been done for the DG method solving problem \rf{1.1}. The first mathematical analysis of this original DG method was given by
Lesaint and Raviart \cite{Lesaint}. They showed that the DG scheme
can be solved in an explicit fashion if $\beta$ is constant and \textcolor[rgb]{0.00,0.00,0.00}{the convergence rate is
of $O(h^k)$-order when the $k$th-order polynomials are used.} Later on,
Johnson and Pitkaranta \cite{Johnson} improved this convergence
rate to $O(h^{k+\frac{1}{2}})$-order. Peterson in \cite{Peterson} further proved that \textcolor[rgb]{0.00,0.00,0.00}{under quasi-uniform triangulation condition, the $
O(h^{k+\frac{1}{2}})$-order convergence is sharp}, namely, the \textcolor[rgb]{0.00,0.00,0.00}{standard convergence rate} of the DG method for problem \rf{1.1} is of $O(h^{k+\frac{1}{2}})$-order, also see Richter's work \cite{Richter}.

On the other hand, in a diametrically opposed effort, some improved
error estimates or \textcolor[rgb]{0.00,0.00,0.00}{optimal convergence rate} are achieved on special meshes. In 1988, Richter
\cite{Richter1} showed that, in the two-dimensional case, the
\textcolor[rgb]{0.00,0.00,0.00}{$L_2$-error estimate of} $\|u-u_h\|\leq Ch^{k+1}\|u\|_{k+2}$ can be obtained for semi-uniform triangle meshes with the assumption that all element edges are bounded away from the characteristic direction $\beta$, that
is, the triangulation $T_h$ satisfies $|\beta\cdot {\bf n}(x)|\geq c_0>0$ on all element edges. Cockburn et. al. in
\cite{Cockburn} proved that if $\beta$ is
a constant vector and the triangulation $T_h$ satisfies the so
called {\em flow condition}, then it holds $\|u-u_h\|\leq Ch^{k+1}|u|_{k+1}$.
Furthermore, Cockburn et. al. in \cite{Cockburn1} relaxed the constant vector $\beta$ to variable $\beta(x)$ (also see \cite{Zhang2}) by adding extra regularity condition of $u\in W^{k+1,\infty}(\Omega)$.

Recently, the WG method has attracted much attention in the
field of numerical partial differential equations
\cite{Xie,Gao,Lin1,Lijian,Lin2,Lin3,Wang,Wang2,Wang3,Zhai,Zhang0,Zhang1,Zhang4}. This method was
introduced and analyzed originally by Wang and Ye in \cite{Wang} for second order elliptic
problems. In general, a WG can be considered as an extension of the standard
finite element method or DG method where the classical derivatives are replaced in
the variational equation by the weakly defined derivatives on
discontinuous weak functions. The main feature of this method is that: (1) the
weak finite element function $u_h=\{u_h^0,\,u_h^b\}$ is used \textcolor[rgb]{0.00,0.00,0.00}{where $u_h^0$ is the value of $u_h$ in the interior of element and $u_h^b$ is the value of $u_h$ on the element boundary. Function $u_h^0$ is totally discontinuous on the partition and $u_h^b$ may be independent with $u_h^0$;} (2)
the weak derivatives are introduced as \textcolor[rgb]{0.00,0.00,0.00}{the distributional derivative of the finite element functions}; The readers are referred to articles \cite{Lin1,Lin2,Wang2}
for more detailed explanation of this method and its relation with
other finite element methods. Although, WG methods have been analyzed for various partial differential equations, for example, the elliptic, parabolic, biharmonic, Stokes and Navier-Stokes equations, et.al., to authors' best knowledge, no WG method is presented for first order hyperbolic problem \rf{1.1} in existing literatures.

In this paper, we will present and analyze a new WG method for solving problem \rf{1.1} \textcolor[rgb]{0.00,0.00,0.00}{using the $k$th-order polynomials. We first construct the WG scheme. Then, we study the stability and approximation properties of the WG solution.} For shape regular meshes consisting of arbitrary polygon/polyhedra, we derive the \textcolor[rgb]{0.00,0.00,0.00}{$L_2$-error estimate} for the WG solution $u_h=\{u_h^0,u_h^b\}$:
\begin{equation}
\|u-u^0_h\|\leq Ch^{k+\frac{1}{2}}|u|_{k+1},\;k\geq 0,\la{1.2}
\end{equation}
and for a special class of meshes (see condition \rf{3.10}), we derive the \textcolor[rgb]{0.00,0.00,0.00}{optimal error estimate:}
\begin{equation}
\|u-u^0_h\|\leq Ch^{k+1}|u|_{k+1},\;k\geq 0.\la{1.3}
\end{equation}
Next, we consider the derivative approximation of the WG solution. By establishing a very simple derivative recovery formula for the convection \textcolor[rgb]{0.00,0.00,0.00}{directional derivative} $\partial_\beta u=\beta\cdot\nabla u$, we obtain the superconvergence results:
\begin{equation}
\|\partial_\beta u-R_h(\partial_\beta u)\|\leq Ch^{k+\frac{1}{2}}|u|_{k+1},\,k\geq 0,\la{1.4}
\end{equation}
and for meshes with condition \rf{3.10},
\begin{equation}
\|\partial_\beta u-R_h(\partial_\beta u)\|\leq Ch^{k+1}|u|_{k+1},\,k\geq 0,\la{1.5}
\end{equation}
where $R_h(\partial_\beta u)$ is the recovery values of $\partial_\beta u$, and $R_h(\partial_\beta u)$ can be computed by means of $u_h$ explicitly, element by element. Our work provides an approach to develop the WG method for first order hyperbolic problems.

  \def\msev{\begin{picture}(   200.00000     ,   200.00000     )(   0.0000000     ,   0.0000000  )
    \def\la{\circle*{0.3}} \multiput(   0.00,   0.00)(   0.250,   0.000){400}{\la}
     \multiput( 100.00,   0.00)(   0.000,   0.250){160}{\la}
     \multiput( 100.00,  40.00)(   0.000,   0.250){ 17}{\la}
     \multiput(  66.67,  44.44)(   0.250,   0.000){133}{\la}
     \multiput(  66.67,  44.44)(   0.000,   0.250){266}{\la}
     \multiput(   0.00,  88.89)(   0.237,   0.079){281}{\la}
     \multiput(   0.00,  75.00)(   0.000,   0.250){ 55}{\la}
     \multiput(   0.00,   0.00)(   0.000,   0.250){300}{\la}
     \multiput( 100.00,   0.00)(   0.250,   0.000){ 44}{\la}
     \multiput( 111.11,   0.00)(   0.250,   0.000){355}{\la}
     \multiput( 200.00,   0.00)(   0.000,   0.250){300}{\la}
     \multiput( 111.11,  88.89)(   0.247,  -0.039){359}{\la}
     \multiput( 111.11,  40.00)(   0.000,   0.250){195}{\la}
     \multiput( 100.00,  40.00)(   0.250,   0.000){ 44}{\la}
     \multiput( 100.00,   0.00)(   0.000,   0.250){160}{\la}
     \multiput( 111.11,  88.89)(   0.247,  -0.039){359}{\la}
     \multiput( 200.00,  75.00)(   0.000,   0.250){ 55}{\la}
     \multiput( 200.00,  88.89)(   0.000,   0.250){444}{\la}
     \multiput( 111.11, 200.00)(   0.250,   0.000){355}{\la}
     \multiput( 111.11, 155.56)(   0.000,   0.250){177}{\la}
     \multiput( 111.11,  88.89)(   0.000,   0.250){266}{\la}
     \multiput(   0.00,  88.89)(   0.237,   0.079){281}{\la}
     \multiput(  66.67, 111.11)(   0.000,   0.250){177}{\la}
     \multiput(  66.67, 155.56)(   0.250,   0.000){177}{\la}
     \multiput( 111.11, 155.56)(   0.000,   0.250){177}{\la}
     \multiput( 100.00, 200.00)(   0.250,   0.000){ 44}{\la}
     \multiput(   0.00, 200.00)(   0.250,   0.000){400}{\la}
     \multiput(   0.00,  88.89)(   0.000,   0.250){444}{\la}
     \multiput( 100.00,  40.00)(   0.250,   0.000){ 44}{\la}
     \multiput( 111.11,  40.00)(   0.000,   0.250){195}{\la}
     \multiput( 111.11,  88.89)(   0.000,   0.250){266}{\la}
     \multiput(  66.67, 155.56)(   0.250,   0.000){177}{\la}
     \multiput(  66.67, 111.11)(   0.000,   0.250){177}{\la}
     \multiput(  66.67,  44.44)(   0.000,   0.250){266}{\la}
     \multiput(  66.67,  44.44)(   0.250,   0.000){133}{\la}
     \multiput( 100.00,  40.00)(   0.000,   0.250){ 17}{\la}\end{picture}  }

  \def\meig{\begin{picture}(   200.00000     ,   200.00000     )(   0.0000000     ,   0.0000000  )
    \def\la{\circle*{0.3}} \multiput(   0.00,   0.00)(   0.250,   0.000){400}{\la}
     \multiput(   0.00,   0.00)(   0.250,   0.000){400}{\la}
     \multiput( 100.00,   0.00)(   0.000,   0.250){160}{\la}
     \multiput(  44.44,  66.67)(   0.225,  -0.108){246}{\la}
     \multiput(   0.00,  88.89)(   0.224,  -0.112){198}{\la}
     \multiput(   0.00,  75.00)(   0.000,   0.250){ 55}{\la}
     \multiput(   0.00,   0.00)(   0.000,   0.250){300}{\la}
     \multiput( 100.00,   0.00)(   0.250,   0.000){ 44}{\la}
     \multiput( 111.11,   0.00)(   0.250,   0.000){355}{\la}
     \multiput( 200.00,   0.00)(   0.000,   0.250){300}{\la}
     \multiput( 133.33,  80.00)(   0.249,  -0.019){267}{\la}
     \multiput( 100.00,  40.00)(   0.160,   0.192){208}{\la}
     \multiput( 100.00,   0.00)(   0.000,   0.250){160}{\la}
     \multiput(  88.89, 133.33)(   0.160,  -0.192){277}{\la}
     \multiput( 133.33,  80.00)(   0.249,  -0.019){267}{\la}
     \multiput( 200.00,  75.00)(   0.000,   0.250){ 55}{\la}
     \multiput( 200.00,  88.89)(   0.000,   0.250){444}{\la}
     \multiput( 111.11, 200.00)(   0.250,   0.000){355}{\la}
     \multiput(  88.89, 133.33)(   0.079,   0.237){281}{\la}
     \multiput(   0.00,  88.89)(   0.224,  -0.112){198}{\la}
     \multiput(  44.44,  66.67)(   0.139,   0.208){320}{\la}
     \multiput(  88.89, 133.33)(   0.079,   0.237){281}{\la}
     \multiput( 100.00, 200.00)(   0.250,   0.000){ 44}{\la}
     \multiput(   0.00, 200.00)(   0.250,   0.000){400}{\la}
     \multiput(   0.00,  88.89)(   0.000,   0.250){444}{\la}
     \multiput(  44.44,  66.67)(   0.225,  -0.108){246}{\la}
     \multiput( 100.00,  40.00)(   0.160,   0.192){208}{\la}
     \multiput(  88.89, 133.33)(   0.160,  -0.192){277}{\la}
     \multiput(  44.44,  66.67)(   0.139,   0.208){320}{\la} \end{picture}  }

\begin{figure}[h]
 \begin{center}\begin{picture}( 340,80 )(   0.0000000     ,   0.0000000     )

\put(113,58){$D$}\put(128,58){$E$}\put(12,16){$C$}\put(202,36){$A$}\put(180,10){$B$}
\put(0,0){\setlength{\unitlength}{0.4pt}\meig}
\multiput(85,0)(40,0){2}{\multiput(0,0)(0,40){2}{\setlength{\unitlength}{0.2pt}\meig}}
\put(170,0){\setlength{\unitlength}{0.4pt}\msev}
\multiput(255,0)(40,0){2}{\multiput(0,0)(0,40){2}{\setlength{\unitlength}{0.2pt}\msev}}
 \end{picture}
\end{center}
\centerline{{\bf Fig. 1}. {\footnotesize Four non-compatible, arbitrary-shaped polygonal grids used in the computation of this paper.}}

\end{figure}

We would emphasize that the WG method can use  non-compatible arbitrary-shaped
    polygonal grids, with only maximal size restriction of $h$,
see Fig. 1.
In these grids, an polygon can have some arbitrarily short edges ($A$ in Fig. 1), can be non-convex ($B$ in Fig. 1),
   can have 180 degree internal angles ($C$ in Fig. 1),
   and can have non-common edges when intersecting neighboring polygons ($D$ and $E$ in Fig. 1). Such grids are used in our numerical experiment in Section 5.

The rest of this paper is organized as follows. In Section 2, we construct the WG scheme for problem \rf{1.1} and show the well-posedness \textcolor[rgb]{0.00,0.00,0.00}{of} this WG scheme. In Section 3, we give the \textcolor[rgb]{0.00,0.00,0.00}{supoptimal and optimal error estimates} for the WG approximation under different mesh conditions. Section 4 is devoted to the derivative approximation of the exact solution by using the derivative recovery formula and the corresponding error estimate is given. Numerical examples are provided in Section 5 to support our theoretical analysis. Some conclusions are given in Section 6.

Throughout this paper, we adopt the notations $H^m(D)$ to
indicate the usual Sobolev spaces on subdomain $D\subset \Omega$
equipped with the norm $\|\cdot\|_{m,D}$ and semi-norm
$|\cdot|_{m,D}$, and when $D=\Omega$, we omit the
subscript $D$. The inner product and norm in space $H^0(\Omega)=L_2(\Omega)$ are
denoted by $(\cdot,\cdot)$ and $\|\cdot\|$, respectively. We
use letter $C$ to represent a generic positive constant,
independent of the mesh size $h$.
\section{Weak Galerkin finite element scheme}
\setcounter{section}{2}\setcounter{equation}{0}
Consider problem \rf{1.1}. For well-posedness (see, e.g., \cite{Friedrichs}), we assume that $\beta\in
[W^1_\infty(\Omega)]^d,\,\alpha \in L_\infty(\Omega),\,f\in
L_2(\Omega),\,g\in L_2(\partial\Omega_-)$, and
\begin{eqnarray}
\alpha+\frac{1}{2}\hbox{div}\beta=\sigma(x)\geq
\sigma_0>0,\;\,x\in \Omega\,\label{2.1}.
\end{eqnarray}

Let $T_h=\bigcup\lbrace K \rbrace$ be a partition of
domain $\Omega$ that consists of arbitrary polygons/polyhedra, where the mesh size $h=\max \,h_K$, $h_K$ is the diameter of element $K$. Assume that the partition $T_h$ is shape regular defined by a set of conditions given in \cite{Wang2}.

First, let us recall the concepts of weak function and weak divergence (see, e.g.,\cite{Wang,Wang2}) which will then be employed to derive a weak \textcolor[rgb]{0.00,0.00,0.00}{Galerkin} finite element scheme for problem \rf{1.1}.
A weak function on element $K$ refers to a function $v=\{ v^0,v^b\}$ with $v^0=v|_{K}\in
L_2(K)$ and $v^b=v|_{\partial K}\in L_2(\partial K)$. Note that for
a weak function $v=\{v^0,v^b\}$, $v^b$ may not be necessarily the trace of $v^0$ on element boundary $\partial K$.

Introduce the weak Galerkin finite element spaces on partition $T_h$:
\begin{eqnarray*}
&&V_h=\{v=\{v^0,v^b\}:\;v^0|_K\in P_k(K),\,v^b|_e\in P_k(e),\,e\subset\partial K,\; K\in
T_h\},\;k\geq 0,\\
&&V_h^0=\{v=\{v^0,v^b\}\in V_h:\, v^b|_{\partial\Omega_-}=0\},
\end{eqnarray*}
where $P_k(D)$ is the space composed of all polynomials on a set $D$ with degree no more than
$k$. We emphasize that in space $V_h$, $v^b$ is set to be single valued on each edge/face $e$. On the other hand, the component $v^0$ is defined element-wise and completely discontinuous on $T_h$. In a certain sense, a weak finite element function $v=\{v^0,v^b\}\in V_h$ is
formed with its components inside all elements glued together by its components on all edges.

For given $v=\{v^0,v^b\}\in V_h$, we define a {\em weak divergence} $\nabla_w\cdot(\beta v)|_K\in P_k(K)$ related to $\beta$ on $K\in T_h$ as the unique solution of the following equation:
\begin{equation}
\int_K\nabla_w\cdot(\beta v)q dx=-\int_Kv^0\beta\cdot\nabla qdx+\int_{\partial K}\beta\cdot{\bf n}v^bqds,\;\forall\,q\in P_k(K),\;K\in T_h.\la{2.2}
\end{equation}

Given a function $u$ with sufficient regularity, one can find an approximation to $u$ by either interpolation or projection in a standard finite element space. In the WG space $V_h$, we will use a locally defined projection of $u$ as its basic approximation. Specifically, let $Q_0:u\in L_2(K)\rightarrow Q_0 u\in P_k(K)$
be the local $L_2$ projection operator such that
\begin{equation}
(u-Q_0u,q)_{K}=0,\;\forall\,q\in P_k(K),\,K\in T_h.\la{2.3}
\end{equation}
It can be shown that \textcolor[rgb]{0.00,0.00,0.00}{\cite[Theorem 1.1]{Zhang2}}
\begin{equation}
\|u-Q_0u\|_{s,K}\leq Ch_K^{m-s}\|u\|_{m,K},\;0\leq s\leq m\leq k+1.\la{2.4}
\end{equation}
The $L_2$ projection operator $Q_b: L_2(e)\rightarrow P_{k}(e),\,e\subset \partial K$ can be similarly defined on the edges of element
$K\in T_h$. Now, we define a projection operator $Q_h: u\in H^1(\Omega)\rightarrow Q_hu\in V_h$ by its action on each element $K$ such that
\begin{equation}
Q_hu|_K=\{Q_0u,Q_bu\},\;K\in T_h.\la{2.5}
\end{equation}

Let $\partial\Omega_\pm=\{x\in\partial\Omega:\,\pm\beta(x)\cdot{\bf n}(x)>0\}$ and $\partial K_\pm=\{x\in\partial K:\,\pm\beta(x)\cdot{\bf n}(x)>0\}$. Then, $\partial\Omega_+$ ($\partial\Omega_-$) is called the outflow (inflow) boundary of domain $\Omega$ with respect to vector $\beta$, similarly, $\partial K_+ \,(\partial K_-)$ is called the outflow (inflow) boundary of element $K$. Furthermore, a face $e$ is called the outflow (inflow) face of element $K$ if $e\subset\partial K_+\;(\partial K_-)$. Denote the sets $\partial T_h=\bigcup \{\partial K: K\in T_h\}$ and $\partial_\pm T_h=\bigcup \{\partial K_\pm: K\in T_h\}$.  For simplicity, we use the following notations,
\begin{eqnarray*}
&&(w,v)_{T_h}=\sum_{K\in T_h}(w,v)_K,\;\;\;\langle w,v\rangle_{\partial T_h}=\sum_{K\in T_h}\langle w,v\rangle_{\partial K}=\sum_{K\in T_h}\int_{\partial K}wvds,\\
&&\langle w,v\rangle_{\partial_\pm T_h}=\sum_{K\in T_h}\langle w,v\rangle_{\partial K_\pm}=\sum_{K\in T_h}\int_{\partial K_\pm}wvds,\;\;\langle w,v\rangle_{\partial \Omega_\pm}=\int_{\partial\Omega_\pm}wvds.
\end{eqnarray*}
Introduce the bilinear form for $w=\{w^0,w^b\},v=\{v^0,v^b\}\in V_h$,
\begin{equation}
a(w,v)=(\nabla_w\cdot(\beta w),v^0)_{T_h}+(\alpha w^0,v^0)_{T_h}+s(w,v),\la{2.6}
\end{equation}
where \textcolor[rgb]{.00,0.00,0.00}{the stabilizer}
$$
s(w,v)=\langle \beta\cdot{\bf n}(w^0-w^b),v^0-v^b\rangle_{\partial_+T_h}.
$$
{\sc Weak Galerkin Method}: a weak Galerkin finite element approximation for problem \rf{1.1} is to find $u_h=\{u_h^0,u_h^b\}\in V_h$ with $u_h^b|_{\partial\Omega_-}=Q_bg$ such that
\begin{equation}
a(u_h,v)=(f,v^0)_{T_h},\;\forall\,v=\{v^0,v^b\}\in V^0_h.\la{2.7}
\end{equation}
\textcolor[rgb]{0.00,0.00,.00}{{\bf Remark 2.1}\quad It should be pointed out that when $\beta\cdot {\bf n}=0$ on some edge/face $e$, the WG equation (2.7) does not contain the unknown $u_h^b|_e$. In this case, we will eliminate the degree of freedom $v^b|_e$ in space $V_h$ or the corresponding basis functions on edge $e$ in space $V_h$.}\\
\yl{ 2.1}\quad{\em For $v=\{v^0,v^b\}\in V_h^0$, it holds}
\begin{equation}
a(v,v)=(\sigma v^0,v^0)_{T_h}+\frac{1}{2}\langle|\beta\cdot{\bf n}|(v^0-v^b),v^0-v^b\rangle_{\partial T_h}+\frac{1}{2}\langle|\beta\cdot{\bf n}|v^b,v^b\rangle_{\partial\Omega_+}.\la{2.8}
\end{equation}
\zm\quad By definition \rf{2.2} and the Green's formula, we have for $w,v\in V_h^0$,
\begin{eqnarray*}
&&(\nabla_w\cdot(\beta w),v^0)_K=-(w^0,\beta\cdot\nabla v^0)_K+\langle\beta\cdot{\bf n}w^b,v^0\rangle_{\partial K}\\
&=&(\nabla\cdot\beta w^0,v^0)_K+(\beta\cdot\nabla w^0,v^0)_K-\langle\beta\cdot{\bf n}w^0,v^0\rangle_{\partial K}+\langle\beta\cdot{\bf n}w^b,v^0\rangle_{\partial K}\\
&=&(\nabla\cdot\beta w^0,v^0)_K-(\nabla_w\cdot(\beta v),w^0)_K+\langle\beta\cdot {\bf n}v^b,w^0\rangle_{\partial K}-\langle\beta\cdot{\bf n}(w^0-w^b),v^0\rangle_{\partial K}\\
&=&(\nabla\cdot\beta w^0,v^0)_K-(\nabla_w\cdot(\beta v),w^0)_K-\langle\beta\cdot{\bf n}(w^0-w^b),v^0-v^b\rangle_{\partial K}
+\langle\beta\cdot {\bf n}v^b,w^b\rangle_{\partial K}.
\end{eqnarray*}
Which implies
$$
(\nabla_w\cdot(\beta v),v^0)_{T_h}=\frac{1}{2}(\nabla\cdot\beta v^0,v^0)_{T_h}-\frac{1}{2}\langle \beta\cdot{\bf n}(v^0-v^b),v^0-v^b\rangle_{\partial T_h}+\frac{1}{2}\langle \beta\cdot {\bf n}v^b,v^b\rangle_{\partial\Omega_+},
$$
where we have used the fact that for $v\in V^0_h$, $v^b$ is single valued on each face $e\subset\partial K$ and $v^b|_{\partial\Omega_-}=0$ such that $\langle \beta\cdot {\bf n}v^b,v^b\rangle_{\partial T_h}=\langle \beta\cdot {\bf n}v^b,v^b\rangle_{\partial\Omega_+}$. Hence, it yields from \rf{2.6} that
\begin{eqnarray*}
a(v,v)&=&(\sigma v^0,v^0)_{T_h}-\frac{1}{2}\langle \beta\cdot{\bf n}(v^0-v^b),v^0-v^b\rangle_{\partial T_h}+\frac{1}{2}\langle \beta\cdot {\bf n}v^b,v^b\rangle_{\partial\Omega_+}+s(v,v)\\
&=&(\sigma v^0,v^0)_{T_h}-\frac{1}{2}\langle \beta\cdot{\bf n}(v^0-v^b),v^0-v^b\rangle_{\partial_- T_h}\\
&&+\frac{1}{2}\langle \beta\cdot{\bf n}(v^0-v^b),v^0-v^b\rangle_{\partial_+ T_h}+\frac{1}{2}\langle \beta\cdot {\bf n}v^b,v^b\rangle_{\partial\Omega_+}\\
&=&(\sigma v^0,v^0)_{T_h}+\frac{1}{2}\langle|\beta\cdot{\bf n}|(v^0-v^b),v^0-v^b\rangle_{\partial T_h}+\frac{1}{2}\langle |\beta\cdot {\bf n}|v^b,v^b\rangle_{\partial\Omega_+}.
\end{eqnarray*}
The proof is completed.\zb

Introduce the norm notation: $|||v|||=a^{\frac{1}{2}}(v,v)$, that is,
\begin{eqnarray}
|||v|||^2=(\sigma v^0,v^0)_{T_h}+\frac{1}{2}\langle|\beta\cdot{\bf n}|(v^0-v^b),v^0-v^b\rangle_{\partial T_h}+\frac{1}{2}\langle |\beta\cdot {\bf n}|v^b,v^b\rangle_{\partial\Omega_+}.\;\la{2.8a}
\end{eqnarray}
\yl{ 2.2}\quad{\em For $v\in V_h^0$, $|||v|||$ defines a norm on space $V_h^0$ and weak Galerkin finite element equation \rf{2.7} has one unique solution.}\\
\zm\quad Let $|||v|||=0$. Then, we have from \rf{2.8a} that $v^0=0$ and $\langle|\beta\cdot{\bf n}|v^b,v^b\rangle_{\partial T_h}=0$ which implies $v^b=0$, therefore, $v=\{v^0,v^b\}=0$ holds such that $|||v|||$ is a norm on space $V_h^0$. Next, let $u_h$ be a solution of problem \rf{2.7} with $f=g=0$. Then, $u_h\in V_h^0$ and $a(u_h,u_h)=|||u_h|||^2=0$, that is, the linear homogeneous problem \rf{2.7} only has zero solution.\zb\\

By Lemma 2.1 and the Cauchy inequality, we can obtain the stability estimate for WG equation \rf{2.7}
$$
|||u_h|||\leq \frac{1}{\sqrt{\sigma_0}}\|f\|\,.
$$

\commentout{
Below, we further examine the structure of equation \rf{2.7} so that we can derive a local solution method. By \rf{2.2}, we have for $v\in V_h^0$,
$$
(\nabla_w\cdot(\beta u_h),v^0)_{T_h}=-(u_h^0,\beta\cdot\nabla v^0)_{T_h}+\langle\beta\cdot{\bf n}u^b_h,v^0\rangle_{\partial T_h},
$$
and
\begin{eqnarray*}
&&\langle\beta\cdot{\bf n}u^b_h,v^0\rangle_{\partial T_h}+s(u_h,v)\\
&=&\langle\beta\cdot{\bf n}u_h^b,v^0-v^b\rangle_{\partial T_h}+\langle\beta\cdot{\bf n}u_h^b,v^b\rangle_{\partial\Omega_+}+\langle\beta\cdot{\bf n}(u_h^0-u^b_h),v^0-v^b\rangle_{\partial_+T_h}\\
&=&\langle\beta\cdot{\bf n}u_h^b,v^0-v^b\rangle_{\partial_-T_h}+\langle\beta\cdot{\bf n}u_h^0,v^0-v^b\rangle_{\partial_+T_h}+\langle\beta\cdot{\bf n}u_h^b,v^b\rangle_{\partial\Omega_+}.
\end{eqnarray*}
Substituting these two equalities into equation \rf{2.7}, we can write WG equation \rf{2.7} in the following form
\begin{eqnarray}
&&-(u_h^0,\beta\cdot\nabla v^0)_{T_h}+(\alpha  u_h^0,v^0)_{T_h}+\langle\beta\cdot{\bf n}u_h^b,v^0-v^b\rangle_{\partial_-T_h}
\nonumber\\
&+&\langle\beta\cdot{\bf n}u_h^0,v^0-v^b\rangle_{\partial_+T_h}+\langle\beta\cdot{\bf n}u_h^b,v^b\rangle_{\partial\Omega_+}=(f,v^0)_{T_h},\;\forall\,v=\{v^0,v^b\}\in V_h^0.\la{2.9}
\end{eqnarray}
Let $\partial K_0=\{x\in\partial K:\,\beta(x)\cdot{\bf n}(x)=0\}$. We assume that partition $T_h$ is such that if an edge/face $e\not\subset \partial K_0$, there must be either $e\subset\partial K_-$ or $e\subset\partial K_+$. Obviously, if  $\beta$ is a constant vector, such assumption always holds for any partition $T_h$. Note that if an edge/face  $e\subset\partial K_+\setminus\partial\Omega_+$, then there must exist an adjacent element $K'$ such that $e\subset\partial K'_-$, see Fig. 1. In \rf{2.9}, taking $v^b=0$ and $v^0(x)=0$ if $x\notin K$, we obtain
\begin{eqnarray}
&&-(u_h^0,\beta\cdot\nabla v^0)_{K}+(\alpha  u_h^0,v^0)_{K}
+\langle\beta\cdot{\bf n}u_h^0,v^0\rangle_{\partial K_+}\nonumber\\
&=&(f,v^0)_{K}-\langle\beta\cdot{\bf n}u^b_h,v^0\rangle_{\partial K_-},\;v^0\in P_k(K).\la{2.10}
\end{eqnarray}
Let $e$ be an edge/face lying on $\partial K_+$ and $v^b(x)=0$ for $x\notin e$. If $e=\partial K_+\cap\partial K'_-$ is an interior edge (noting that ${\bf n'}|_{e\cap \partial K'_-}=-{\bf n}|_{e\cap\partial K_+}$), it holds
$$
-\langle\beta\cdot{\bf n}u_h^b,v^b\rangle_{\partial_-T_h}=\langle\beta\cdot{\bf n}u_h^b,v^b\rangle_{\partial K_+\cap e},\,e\not\subset\partial\Omega_+,
$$
otherwise $\langle\beta\cdot{\bf n}u_h^b,v^b\rangle_{\partial_-T_h}=0,\,e\subset\partial\Omega_+$. Therefore, taking $v^0=0$ and $v^b(x)=0$ if $x\notin e$, we have from \rf{2.9} for $e\subset \partial K_+$,
\begin{eqnarray}
\langle\beta\cdot{\bf n}u_h^b,v^b\rangle_{\partial K_+\cap e}
 =\langle\beta\cdot{\bf n}u_h^0,v^b\rangle_{\partial K_+\cap e},\;v^b\in P_k(e).\ \ \ \la{2.11}
\end{eqnarray}
Equations \rf{2.10}--\rf{2.11} are the local form of weak Galerkin finite element equation \rf{2.7} restricted on element $K$.

Then, we can write WG equation \rf{2.7} in the following form
\begin{eqnarray}
&&-(u_h^0,\beta\cdot\nabla v^0)_{T_h}+(\alpha  u_h^0,v^0)_{T_h}+\langle\beta\cdot{\bf n}u_h^b,v^0\rangle_{\partial T_h}
\nonumber\\
&+&\langle\beta\cdot{\bf n}(u_h^0-u_h^b),v^0-v^b\rangle_{\partial_+T_h}=(f,v^0)_{T_h},\;\forall\,v=\{v^0,v^b\}\in V_h^0.\la{2.9}
\end{eqnarray}
In \rf{2.9}, we first take $v^b=0$ and $v^0(x)=0$ if $x\notin K$, then it yields
\begin{eqnarray}
&&-(u_h^0,\beta\cdot\nabla v^0)_{K}+(\alpha  u_h^0,v^0)_{K}
+\langle\beta\cdot{\bf n}u_h^0,v^0\rangle_{\partial K_+}\nonumber\\
&=&(f,v^0)_{K}-\langle\beta\cdot{\bf n}u^b_h,v^0\rangle_{\partial K_-},\;v^0\in P_k(K).\la{2.10}
\end{eqnarray}
Next, taking $v^0=0$ and $v^b(x)=0$ if $x\notin e$ in \rf{2.9}, we obtain
\begin{eqnarray}
\langle\beta\cdot{\bf n}u_h^b,v^b\rangle_{e}
 =\langle\beta\cdot{\bf n}u_h^0,v^b\rangle_{e},\;v^b\in P_k(e),\,e\subset\partial_+T_h.\ \ \ \la{2.11}
\end{eqnarray}
Equations \rf{2.10} and \rf{2.11} are the local form of weak Galerkin finite element equation \rf{2.7} restricted on element $K$.\\
{\bf Remark 2.2} Equation \rf{2.11} does not imply $u_h^b|_e=u_h^0|_{e}$ like the original DG method. For example, let $e\subset\partial K'_-\cap\partial_+T_h$ be an interior inflow face of element $K'$, noting that $e$ also may lie on the outflow faces of some other elements $\{K\}$ where $e\cap\partial K_+\not=\varnothing$. Then, when $e$ contains a hanging point such that it is composed of two outflow faces of elements in $\{K\}$, $u_h^0|_e=u_h^0|_{e\cap\partial_+ T_h}$ is composed of two different polynomials, but $u_h^b|_e$ always is a signal polynomial.
\begin{center}
\scalebox{0.5}{\includegraphics{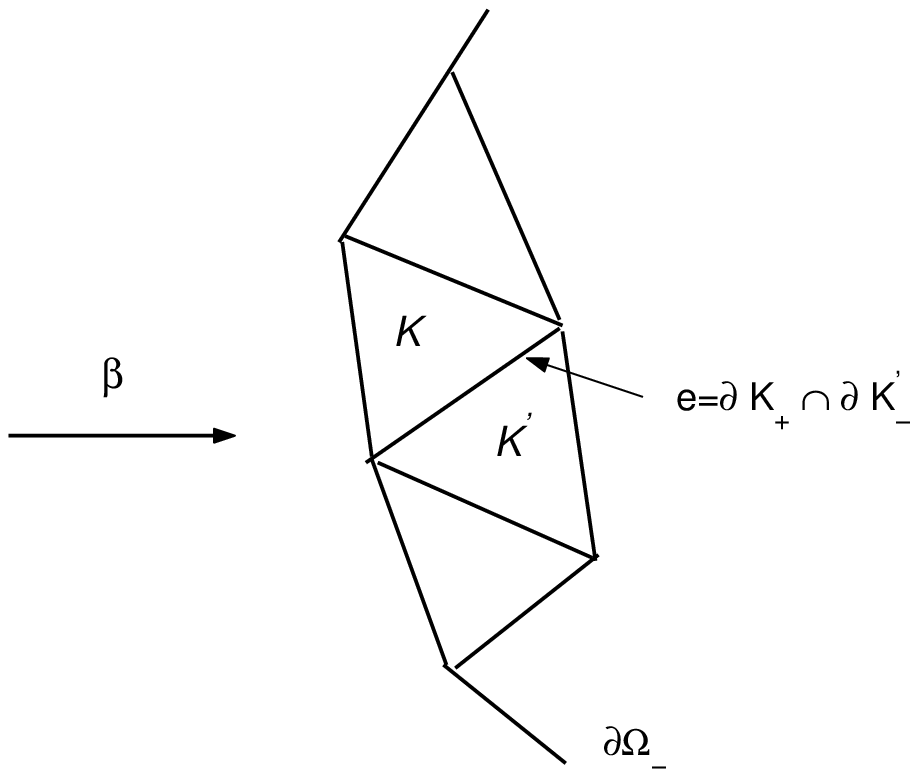}}
\centerline{\small {\bf
Fig.2}\quad Inflow element $K$ and its adjacent element $K'$}
\end{center}
Now, let us show how discrete equation \rf{2.7} can be solved locally, element by element. We call $K$ an inflow  element if $\partial K_-\subset\partial\Omega_-$, see Fig. 2. We first solve those values of $\{u_h^0\}$ on the inflow elements $\{K\}$ by equation \rf{2.10}, element by element, noting that $u_h^b|_{\partial K_-}=Q_bg$ is known. Then, from equation \rf{2.11}, we can solve the values of $\{u_h^b\}$ on boundaries $\{\partial K_+\}$ of inflow elements, edge by edge. Thus, we obtain $u_h=\{u_h^0,u_h^b\}$ on the inflow elements $\{K\}$. Next, for those elements $\{K'\}$ which are adjacent to some inflow element $K$ (its an edge $e\cap\partial K'_-\subset\partial K_+$, see Fig.2), we can solve $\{u_h^0\}$ by equation \rf{2.10} on each element $K'$, and then solve $\{u_h^b\}$ on each edge $e\cap\partial K'_+$ by equation \rf{2.11}. Again, we obtain $u_h=\{u_h^0,u_h^b\}$ on elements $\{K'\}$. Then, we can further consider $u_h$ on elements which are adjacent to some element $K'$. In such a sequence, finally we can find $u_h=\{u_h^0,u_h^b\}$ locally, element by element, for all $K\in T_h$.\\

From the analysis above, we see that the WG method has an explicit structure that is very similar to the original DG method for solving problem \rf{1.1}, cf. \cite{Lesaint}. In the appendix, we will give a brief discussion on the relation between the WG method and the original DG method.
}

\section{Error Analysis}
\setcounter{section}{3}\setcounter{equation}{0}
In this section, we do the error analysis for the WG method. We will derive the \textcolor[rgb]{0.00,0.00,.00}{supoptimal $L_2$-error estimates for general mesh and optimal $L_2$-error estimate for special mesh.}
\subsection{Error estimate on general meshes}
For any function $u\in H^1(K)$, the following trace inequality holds.
\begin{equation}
\|u\|^2_{L_2(e)}\leq C\big(h_K^{-1}\|u\|^2_{0,K}+h_K\|\nabla u\|_{0,K}^2\big),\;e\subset\partial K,\,K\in T_h.\la{3.1}
\end{equation}
\yl{ 3.1}\quad{\em Let $u\in H^1(\Omega)$. Then it holds true
\begin{equation}
(\nabla\cdot(\beta u),v^0)=(\nabla_w\cdot(\beta Q_hu),v^0)_{T_h}-l_1(u,v)+l_2(u,v),\,\forall\,v=\{v^0,v^b\}\in V_h^0,\la{3.2}
\end{equation}
where}
\begin{eqnarray*}
&&l_1(u,v)=(u-Q_0u,\beta\cdot\nabla v^0)_{T_h},\\
&&l_2(u,v)=\langle \beta\cdot{\bf n}(u-Q_bu),v^0-v^b\rangle_{\partial T_h}+\langle \beta\cdot{\bf n}(u-Q_bu),v^b\rangle_{\partial\Omega_+}.
\end{eqnarray*}
\zm\quad From the Green's formula, definition \rf{2.5} of projection $Q_h$ and definition \rf{2.2} of the weak divergence, we have
\begin{eqnarray*}
(\nabla\cdot(\beta u),v^0)&=&-(u,\beta\cdot\nabla v^0)_{T_h}+
\langle \beta\cdot{\bf n}u,v^0\rangle_{\partial T_h}\\
&=&-(Q_0u,\beta\cdot\nabla v^0)_{T_h}-l_1(u,v)+
\langle \beta\cdot{\bf n}Q_bu,v^0\rangle_{\partial T_h}-
\langle \beta\cdot{\bf n}Q_bu,v^0\rangle_{\partial T_h}\\
&&+\langle \beta\cdot{\bf n}u,v^0-v^b\rangle_{\partial T_h}+\langle \beta\cdot{\bf n}u,v^b\rangle_{\partial T_h}\\
&=&(\nabla_w\cdot(\beta Q_hu),v^0)_{T_h}-l_1(u,v)+\langle \beta\cdot{\bf n}(u-Q_bu),v^0\rangle_{\partial T_h}\\
&&-\langle \beta\cdot{\bf n}(u-Q_bu),v^b\rangle_{\partial T_h}+\langle \beta\cdot{\bf n}(u-Q_bu),v^b\rangle_{\partial T_h}\\
&=&(\nabla_w\cdot(\beta Q_hu),v^0)_{T_h}-l_1(u,v)+l_2(u,v),
\end{eqnarray*}
where we have used the fact $v^b|_{\partial\Omega_-}=0$. \zb

Using Lemma 3.1, equation \rf{1.1} and definition \rf{2.6} of $a(u,v)$, we immediately obtain the following result.\\
\yl{ 3.2}\quad{\em Let $u\in H^1(\Omega)$ be the solution of problem \rf{1.1}. Then we have
\begin{equation}
a(Q_hu,v)=(f,v^0)_{T_h}+l_1(u,v)-l_2(u,v)+l_3(u,v)+s(Q_hu,v),\;\forall\,v\in V_h^0,\la{3.3}
\end{equation}
where $l_3(u,v)=(\alpha(Q_0u-u),v^0)_{T_h}$}.\\
\textcolor[rgb]{0.00,0.00,0.00}{\zm\quad From equation \rf{1.1} and Lemma 3.1, we obtain
\begin{eqnarray*}
(\nabla_w\cdot(\beta Q_hu),v^0)_{T_h}+(\alpha u,v^0)_{T_h}-l_1(u,v)+l_2(u,v)=(f,v^0)_{T_h},
\end{eqnarray*}
or
\begin{eqnarray*}
&&(\nabla_w\cdot(\beta Q_hu),v^0)_{T_h}+(\alpha Q_0u,v^0)_{T_h}+s(Q_hu,v)\\
&=&(f,v^0)_{T_h}+l_1(u,v)-l_2(u,v)+(\alpha(Q_0u-u),v^0)_{T_h}+s(Q_hu,v).
\end{eqnarray*}
Combining this with the definition \rf{2.6} of $a(u,v)$, the proof is completed.\zb}

Now, we are in the position to derive the error estimate for the WG solution $u_h$.\\
\dl{ 3.1}\quad{\em Assume that $T_h$ is a shape regular partition and let $u\in H^{k+1}(\Omega)$ ($k\geq 0$) and $u_h\in V_h$ be the solutions of problems \rf{1.1} and the WG equation \rf{2.7}, respectively. Then, the following error estimate holds.}
\begin{equation}
|||Q_hu-u_h|||\leq Ch^{k+\frac{1}{2}}|u|_{k+1},\;k\geq 0.\la{3.4}
\end{equation}
\zm\quad \textcolor[rgb]{.00,0.00,0.00}{Denote the error function by $e_h=Q_hu-u_h\in V_h^0$.} From equation \rf{2.7} and \rf{3.3}, we obtain
\begin{eqnarray}
a(e_h,v)=l_1(u,v)-l_2(u,v)+l_3(u,v)+s(Q_hu,v),\;v\in V_h^0.\la{3.5}
\end{eqnarray}
Below we estimate terms $l_i(u,v)$ ($i=1,2,3$) and $s(Q_h,v)$. It follows from the definition of $Q_0$ and \textcolor[rgb]{.00,0.00,0.00}{the finite element inverse inequality: $\|\nabla v^0\|_{0,K}\leq Ch^{-1}_K\|v^0\|_{0,K}$,
\begin{eqnarray}
&&|l_1(u,v)+l_3(u,v)|=|(u-Q_0u,(\beta-\beta^c)\cdot\nabla v^0)_{T_h}+(\alpha(u-Q_0u),v^0)_{T_h}|\nonumber\\
&\leq& C\sum_{K\in T_h}h_K^{k+1}|u|_{k+1,K}h_K\|\nabla v^0\|_{0,K}+Ch^{k+1}|u|_{k+1}\|v^0\|
\leq Ch^{k+1}|u|_{k+1}|||v|||,\ \ \ \ \ \ \la{3.6}
\end{eqnarray}
}where $\beta^c$ is the piecewise constant approximation of $\beta$ on $T_h$. Next, it follows from the definition of $Q_b$,
\begin{equation}
\|u-Q_bu\|_{L_2(e)}^2=\int_e(u-Q_bu)(u-Q_0u)ds\leq \|u-Q_bu\|_{L_2(e)}\|u-Q_0u\|_{L_2(e)},\,e\subset\partial K,\nonumber
\end{equation}
which implies $\|u-Q_bu\|_{L_2(e)}\leq \|u-Q_0u\|_{L_2(e)}$. Then, we have by using the trace inequality,
\begin{eqnarray}
|l_2(u,v)|&=&\langle \beta\cdot{\bf n}(u-Q_bu),v^0-v^b\rangle_{\partial T_h}+\langle \beta\cdot{\bf n}(u-Q_bu),v^b\rangle_{\partial\Omega_+}\nonumber\\
&\leq& \Big(|\beta|_\infty\sum_{K\in T_h}\|u-Q_0u\|^2_{0,\partial K}\Big)^{\frac{1}{2}}|||v|||\leq Ch^{k+\frac{1}{2}}|u|_{k+1}|||v|||,\la{3.7}
\end{eqnarray}
and
\begin{eqnarray}
&&s(Q_hu,v)=\langle\beta\cdot{\bf n}(Q_0u-u+u-Q_bu),v^0-v^b\rangle_{\partial_+T_h}\nonumber\\
&\leq& \Big(\sum_{K\in T_h}\|\,|\beta\cdot{\bf n}|^{\frac{1}{2}}(Q_0u-u+u-Q_bu)\|_{0,\partial K}^2\Big)^{\frac{1}{2}}
\Big(\sum_{K\in T_h}\|\,|\beta\cdot{\bf n}|^{\frac{1}{2}}(v^0-v^b)\|_{0,\partial K}^2\Big)^{\frac{1}{2}}\nonumber\\
&\leq& \Big(\sum_{K\in T_h}2|\beta|_\infty\|Q_0u-u\|_{0,\partial K}^2\Big)^{\frac{1}{2}}\,|||v|||\leq Ch^{k+\frac{1}{2}}|u|_{k+1}|||v|||.\la{3.8}
\end{eqnarray}
Substituting \rf{3.6}--\rf{3.8} into \rf{3.5}, it yields
$$
|a(e_h,v)|\leq Ch^{k+\frac{1}{2}}|u|_{k+1}|||v|||,\;v\in V_h^0.
$$
Taking $v=e_h\in V_h^0$ and using Lemma 2.1, it yields
$$
|||e_h|||^2=a(e_h,e_h)\leq Ch^{k+\frac{1}{2}}|u|_{k+1}|||e_h|||.
$$
The proof is completed.\zb

From Theorem 3.1 and the triangle inequality, we immediately obtain the \textcolor[rgb]{0.00,0.00,.00}{supoptimal}  $L_2$-error estimate,
\begin{equation}
\|u-u_h^0\|\leq Ch^{k+\frac{1}{2}}|u|_{k+1},\;k\geq 0.\la{3.9}
\end{equation}
\subsection{Error estimate on special meshes}
In this subsection, we consider the \textcolor[rgb]{0.00,0.00,.00}{optimal $L_2$-error estimate} of the WG solution under the special mesh condition.

Denote by $\mathcal{E}_h$ the union of all element edges in $T_h$ and set its subset
\begin{eqnarray*}
\mathcal{E}^0_h=\{e\in \mathcal{E}_h: \hbox{there exists point $x_0\in e\subset\partial K$ such that $|\beta\cdot{\bf n}(x_0)|\leq h_K$}\}
\end{eqnarray*}
Obviously, if $x_0\in e$ is such that $\beta(x_0)=0$ or vector $\beta(x_0)$ is (almost) parallel to edge $e$, then $e\in\mathcal{E}_h^0$ holds true.
\commentout{\yl{ 3.3}\quad{\em Let $e\subset\partial K$ be an edge of element $K$ and assume that there exists a point $x_0\in e$ such that $\beta\cdot{\bf n}(x_0)>0$. Then, $e=e_K^+$ is an outflow face of element $K$, otherwise edge $e$ is in  $\mathcal{E}^0_h$.}\\
\zm\quad Let point $x_0\in e$ such that $\beta\cdot{\bf n}(x_0)>0$. If $\beta\cdot{\bf n}(x)>0,\;\forall\,x\in e$, then $e=e_K^+$ is an outflow face. Otherwise, there is a point $x_1\in e$ such that $\beta\cdot{\bf n}(x_1)\leq 0$. If $\beta\cdot{\bf n}(x_1)=0$, then $e$ is in $\mathcal{E}^0_h$. If $\beta\cdot{\bf n}(x_1)<0$, according to the intermediate value theorem of continuous function, there must exists a point $x_2\in e$ between $x_0$ and $x_1$ such that $\beta\cdot{\bf n}(x_2)=0$. Therefore, edge $e$ is in $\mathcal{E}^0_h$.\zb
}
\textcolor[rgb]{0.00,0.00,.00}{It is possible that a face $e$ is both an outflow face and a face in $\mathcal{E}^0_h$.} \commentout{From Lemma 3.3, it is easy to see that if set $\partial_+T_h\backslash\mathcal{E}^0_h$ is no empty, then it is composed of some outflow faces of elements in $T_h$.}

Now, we introduce the following mesh condition for partition $T_h$:
\begin{eqnarray}
\hbox{Each element $K$ has at most one outflow face $e_K^+$ in set $\partial K\backslash\mathcal{E}^0_h$}.\;\;\;\;\label{3.10}
\end{eqnarray}
Fig.2 shows a triangulation satisfying the mesh condition \rf{3.10}.
\begin{center}
\scalebox{0.5}{\includegraphics{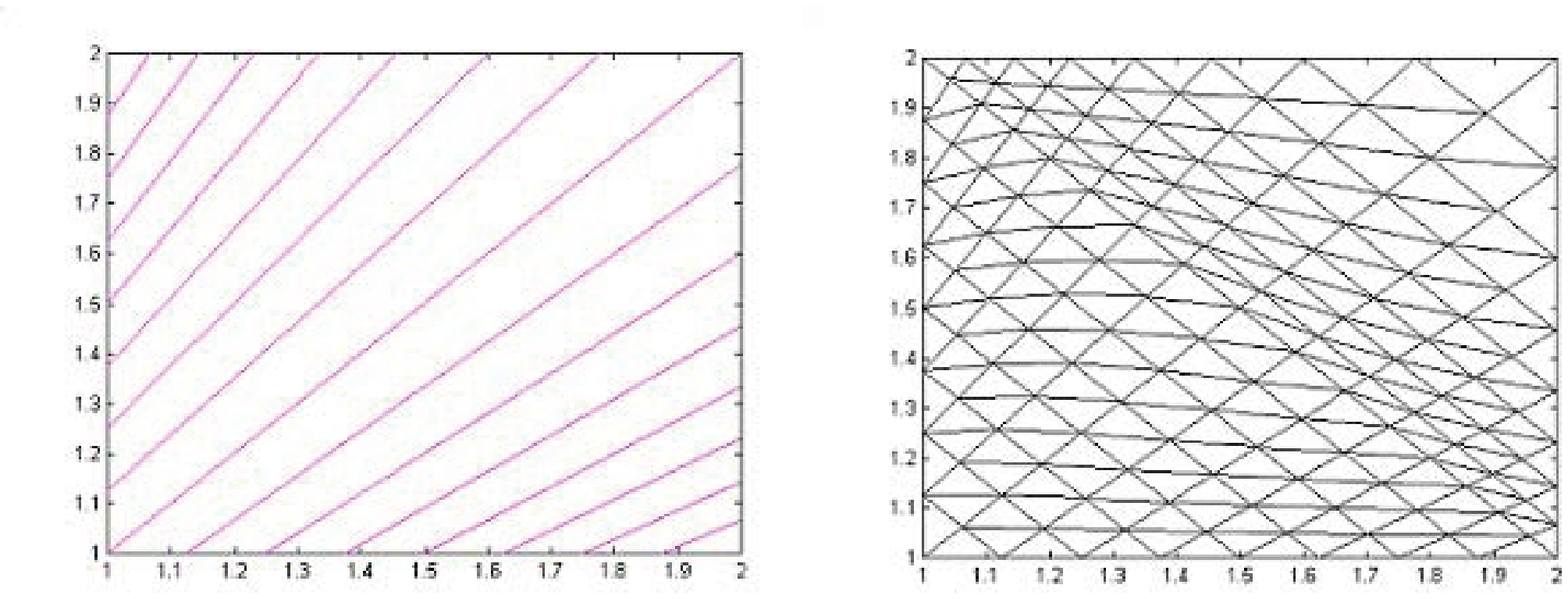}}

{\small {\bf Fig. 2.}\quad
The streamlines of
$\boldsymbol{\beta}=(x,y)$ (left) and the actual mesh satisfying condition \rf{3.10}(right).}
\end{center}

Note that condition \rf{3.10} has no restriction on the number of outflow faces of each element $K$, \textcolor[rgb]{0.00,0.00,.00}{some extra outflow faces (if exist) can be in $\mathcal{E}^0_h$.}\\
{\bf Remark 3.1}\quad A requirement in the mesh condition ({\em flow condition}) given in \cite{Cockburn,Cockburn1} for the DG method is that:  "Each interior outflow face $e^+_K$ is included in an inflow face of another element". This requirement implies that no hanging points are allowed on interior outflow faces. Condition \rf{3.10} has no such requirement and allows the partition with hanging points on outflow/inflow faces.\\
{\bf Remark 3.2}\quad \textcolor[rgb]{0.00,0.00,.00}{Let $\mathcal{E}_h^+=\{e_K\subset\partial K_+: \hbox{element $K$ has at least two outflow faces}\}$ and $\mathcal{E}_h^c=\mathcal{E}_h^+/\mathcal{E}_h^0$.} The {\em flow condition} in \cite{Cockburn1} also contains the following requirement:
\begin{equation}
\sum_{e_K\in\mathcal{E}_h^c}|e_K|^{\frac{d-2}{d-1}}\leq C\,.\la{3.10a}
\end{equation}
When $\mathcal{E}_h^c\not=\emptyset$, using \rf{3.10a} and the additional assumption of $u\in W^{k+1,\infty}(\Omega)$, article \cite{Cockburn1} derived the $(k+1)$th-order convergence. Under conditions of \rf{3.10}-\rf{3.10a} and $u\in W^{k+1,\infty}(\Omega)$, we also can derive the same result. However, it is easy to make a partition \textcolor[rgb]{.00,0.00,0.00}{to satisfy} $\mathcal{E}_h^c=\emptyset$. Therefore,
to show the optimality on both the convergence order and the regularity requirement (see Theorem 3.2), we have not considered the situation of $\mathcal{E}_h^c\not=\emptyset$.

\textcolor[rgb]{0.00,0.00,.00}{To obtain the optimal error estimate}, we still need to
introduce a special projection. For $u\in H^1(\Omega)$, define the projection operator $P^+_h$, retracted on each element $K\in T_h$, $P_h^+u\in P_k(K)$ such that
\begin{eqnarray}
&&\int_K(u-P_h^+u)v\,dx=0,\;\;\forall\, v\in P_{k-1}(K),\label{3.11}\\
&&\int_{e^+_K}(u-P_h^+u)v\,ds=0,\;\;\forall\, v\in
P_k(e^+_K),\label{3.12}
\end{eqnarray}
where condition \rf{3.11} is
vacuous if $k=0$. \textcolor[rgb]{0.00,0.00,.00}{In \rf{3.12}, $e^+_K\subset\partial K$ is the outflow face in set  $\partial _+T_h\backslash\mathcal{E}^0_h$ if such face exists (see condition \rf{3.10}), otherwise we can select some face $e\subset\partial K$ as $e_K^+$ to define the projection $P_h^+u$. This implies that all outflow faces in set $\partial _+T_h\backslash\mathcal{E}^0_h$ have been selected to define $P_h^+u$ in \rf{3.12}.} Projection $P_h^+u$ has been used in some
articles \cite{Cockburn1,Zhang3}. By the argument in \cite{Zhang2,Zhang3}, we have the following result.\\
\textcolor[rgb]{0.00,0.00,.00}{\yl{ 3.3}}\quad{\em The projection function $P_h^+u$ is well posed and satisfies the
approximation property}
\begin{equation}
\|u-P_h^+u\|_{0,K}+h_K^{\frac{1}{2}}\|u-P_h^+u\|_{0,\partial K}\leq
Ch_K^{k+1}|u|_{H^{k+1}(K)},\;k\geq 0,\;\;K\in T_h.\label{3.13}
\end{equation}

By means of $P_h^+$, we can introduce a new projection operator $Q_h^+:\,u\in H^1(\Omega)\rightarrow Q_h^+u\in V_h$ defined by
$$
Q_h^+=\{Q_0^+,Q_b\}\doteq\{P_h^+,Q_b\}.
$$
It is easy to see that all analysis \textcolor[rgb]{0.00,0.00,.00}{maintained above} to hold if $Q_hu$ is replaced by projection $Q_h^+u$.\\
\textcolor[rgb]{0.00,0.00,.00}{\yl{ 3.4}}\quad{\em For $e\subset\partial K$, it holds that when $e\in\mathcal{E}_h^0$,
\begin{eqnarray}
|\beta\cdot{\bf n}(x)|\leq |(|\beta|_{1,\infty}+1)h_K,\;x\in e,\;e\in\mathcal{E}_h^0,\la{3.14}
\end{eqnarray}
\textcolor[rgb]{0.00,0.00,.00}{and when} $e\not\in\mathcal{E}_h^0$,}
\begin{eqnarray}
|\beta\cdot{\bf n}(x)|>h_K,\;x\in e,\,e\not\in\mathcal{E}_h^0.\la{3.15}
\end{eqnarray}
\zm\quad We know that when $e\in \mathcal{E}_h^0$, there exists point $x_0\in e$ such that $|(\beta\cdot{\bf n})(x_0)|\leq h_K$. Since ${\bf n}(x_0)={\bf n}(x)$ on edge $e$, we obtain
\begin{equation}
|\beta\cdot{\bf n}(x)|\leq|(\beta(x)-\beta(x_0))\cdot{\bf n}(x)|+|(\beta\cdot{\bf n})(x_0)|\leq h_K|\beta|_{1,\infty}+h_K,\;x\in e,\nonumber
\end{equation}
this gives \rf{3.14}. Estimate \rf{3.15} comes from the definition of $\mathcal{E}_h^0$.\zb\\
\dl{ 3.2}\quad{\em Assume that $T_h$ is a shape regular partition satisfying the mesh condition \rf{3.10} and let $u\in H^{k+1}(\Omega)$ ($k\geq 0$) and $u_h\in V_h$ be the solutions of problems \rf{1.1} and WG equation \rf{2.7}, respectively. Then, we have the following \textcolor[rgb]{0.00,0.00,.00}{optimal error estimate}.}
\begin{equation}
\|u-u_h^0\|\leq Ch^{k+1}|u|_{k+1},\;k\geq 0.\la{3.16}
\end{equation}
\zm\quad \textcolor[rgb]{0.00,0.00,.00}{Denote the error function by $e^+_h=Q^+_hu-u_h\in V_h^0$.} From equation \rf{3.5}, we obtain
\begin{eqnarray}
a(e^+_h,v)=l_1(u,v)-l_2(u,v)+l_3(u,v)+s(Q^+_hu,v),\;v\in V_h^0,\la{3.17}
\end{eqnarray}
where $l_i(u,v)$ ($i=1,2,3$) are given by \rf{3.2}--\rf{3.3} in which $Q_0$ is replaced by $Q_0^+$. From the \textcolor[rgb]{.00,0.00,0.00}{definition of operator $Q_0^+$ and the finite element inverse inequality}, we have
\begin{eqnarray*}
&&|l_1(u,v)||=|(u-Q^+_0u,(\beta-\beta^c)\cdot\nabla v^0)_{T_h}|\leq Ch^{k+1}|u|_{k+1}|||v|||,\\
&&|l_3(u,v)|=|(\alpha(u-Q_0^+u),v^0)_{T_h}|\leq Ch^{k+1}|u|_{k+1}|||v|||.
\end{eqnarray*}
Next, we write
\begin{eqnarray*}
l_2(u,v)=\langle \beta\cdot{\bf n}(u-Q_bu),v^0-v^b\rangle_{\partial T_h}+\langle \beta\cdot{\bf n}(u-Q_bu),v^b\rangle_{\partial\Omega_+}=S_1(u,v)+S_2(u,v),
\end{eqnarray*}
where
\begin{eqnarray*}
&&S_1(u,v)=\langle \beta\cdot{\bf n}(u-Q_bu),v^0-v^b\rangle_{\mathcal{E}_h^0}+\langle \beta\cdot{\bf n}(u-Q_bu),v^0-v^b\rangle_{\partial T_h\backslash\mathcal{E}_h^0},\\
&&S_2(u,v)=\langle \beta\cdot{\bf n}(u-Q_bu),v^b\rangle_{\partial\Omega_+\bigcap\mathcal{E}_h^0}+\langle \beta\cdot{\bf n}(u-Q_bu),v^b\rangle_{\partial\Omega_+\backslash\mathcal{E}_h^0}.
\end{eqnarray*}
First, it follows from \rf{3.14} that
\begin{eqnarray*}
\langle \beta\cdot{\bf n}(u-Q_bu),v^0-v^b\rangle_{\mathcal{E}_h^0}&\leq& C\sum_{e\in \mathcal{E}_h^0}h_K^{\frac{1}{2}}\|u-Q_bu\|_{L_2(e)}\|\,|\beta\cdot{\bf n}|^{\frac{1}{2}}(v^0-v^b)\|_{L_2(e)}\\
&\leq& Ch^{k+1}|u|_{k+1}|||v|||.
\end{eqnarray*}
Next, it follows from the definition of $Q^b$ and \rf{3.15} that
\begin{eqnarray*}
&&\langle \beta\cdot{\bf n}(u-Q_bu),v^0-v^b\rangle_{\partial T_h\backslash\mathcal{E}_h^0}=\langle (\beta-\beta^c)\cdot{\bf n}(u-Q_bu),v^0-v^b\rangle_{\partial T_h\backslash\mathcal{E}_h^0}\\
&\leq& 2|\beta|_{1,\infty}\sum_{e\in\mathcal{E}_h\backslash\mathcal{E}_h^0}h_K^{\frac{1}{2}}\|u-Q_bu\|_{L_2(e)}
\|h_K^{\frac{1}{2}}(v^0-v^b)\|_{L_2(e)}\\
&\leq& 2|\beta|_{1,\infty}\sum_{e\in\mathcal{E}_h\backslash\mathcal{E}_h^0}h_K^{\frac{1}{2}}\|u-Q_bu\|_{L_2(e)}
\|\,|\beta\cdot{\bf n}|^{\frac{1}{2}}(v^0-v^b)\|_{L_2(e)}\leq Ch^{k+1}|u|_{k+1}|||v|||.
\end{eqnarray*}
Thus, it yields
$$
S_1(u,v)\leq Ch^{k+1}|u|_{k+1}|||v|||.
$$
Similarly, $S_2(u,v)\leq Ch^{k+1}|u|_{k+1}|||v|||$ can be derived. Therefore, we obtain
$$
|l_2(u,v)|\leq Ch^{k+1}|u|_{k+1}|||v|||.
$$
\textcolor[rgb]{0.00,0.00,.00}{Finally, we write $s(Q^+_hu,v)$ as follows},
\begin{eqnarray*}
s(Q^+_hu,v)&=&\langle\beta\cdot{\bf n}(Q_0^+u-Q_bu),v^0-v^b\rangle_{\partial_+T_h}=\langle\beta\cdot{\bf n}(Q_0^+u-Q_bu),v^0-v^b\rangle_{\partial_+T_h\bigcap\mathcal{E}_h^0}\\
&&+\langle\beta\cdot{\bf n}(Q_0^+u-Q_bu),v^0-v^b\rangle_{\partial_+T_h\backslash\mathcal{E}_h^0}=R_1(u,v)+R_2(u,v).
\end{eqnarray*}
First, \textcolor[rgb]{0.00,0.00,.00}{it follows from} \rf{3.13} and \rf{3.14} that
\begin{eqnarray*}
R_1(u,v)&\leq& C\sum_{e\in\mathcal{E}_h^0}h_K^{\frac{1}{2}}\|Q^+_0u-u+u-Q_bu\|_{L_2(e)}\|\,|\beta\cdot{\bf n}|^{\frac{1}{2}}(v^0-v^b)\|_{L_2(e)}\\
&\leq& Ch^{k+1}|u|_{k+1}|||v|||\,.
\end{eqnarray*}
Next, it follows from the definitions of operators $Q^+_0=P_h^+$ and $Q_b$, and \rf{3.15} that
\begin{eqnarray*}
&&R_2(u,v)=\sum_{e_K^+\in\partial_+T_h\backslash\mathcal{E}_h^0}\langle(\beta-\beta^c)\cdot{\bf n}(Q_0^+u-u+u-Q_bu),v^0-v^b\rangle_{e_K^+}\\
&\leq& |\beta|_{1,\infty}\sum_{e\in\mathcal{E}_h\backslash\mathcal{E}_h^0}h_K^{\frac{1}{2}}\|Q^+_0u-u+u-Q_bu\|_{L_2(e)}
\|\,h_K^{\frac{1}{2}}(v^0-v^b)\|_{L_2(e)}\leq Ch^{k+1}|u|_{k+1}|||v|||\,.
\end{eqnarray*}
Together with estimate of $R_1(u,v)$, it yields
$$
s(Q_h^+u,v)\leq Ch^{k+1}|u|_{k+1}|||v|||\,.
$$
Substituting estimates of $l_1(u,v)-l_3(u,v)$ and $s(Q_h^+u,v)$ into \rf{3.17} and taking $v=e_h^+$, we obtain
\begin{equation}
|||e^+_h|||^2=a(e_h^+,e^+_h)\leq Ch^{k+1}|u|_{k+1}|||e^+_h|||,\;k\geq 0,\la{3.18}
\end{equation}
which also implies
\begin{equation}
\|Q^+_0u-u_h^0\|\leq Ch^{k+1}|u|_{k+1},\;k\geq 0.\la{3.19}
\end{equation}
The proof is completed by using the triangle inequality and the approximation property of operator $Q^+_0=P_h^+$.\zb

\section{Recovery formula for the derivative approximation}
\setcounter{section}{4} \setcounter{equation}{0}
In this section, we give the approximation to the convection \textcolor[rgb]{.00,0.00,0.00}{directional derivative} $\partial_\beta u=\beta\cdot\nabla u$ by using a very simple derivative recovery formula.

Let $u_h$ be the WG solution of equation \rf{2.7}. Define the recovery formula of derivative $\partial_\beta u$ by
\begin{equation}
R_h(\partial_\beta u)=f-(\alpha +\nabla\cdot\beta)u_h^0.\la{3.18}
\end{equation}
\dl{ 4.1}\quad{\em Assume that $T_h$ is a shape regular partition and $u\in H^{k+1}(\Omega)$ ($k\geq 0$) is the solution of problem \rf{1.1}. Then, the recovery derivative $R_h(\partial_\beta u)$ given by \rf{3.18} admits the following superconvergence estimate.
\begin{equation}
\|\partial_\beta u-R_h(\partial_\beta u)\|\leq Ch^{k+\frac{1}{2}}|u|_{k+1},\;k\geq 0.\la{3.19}
\end{equation}
Moreover, if partition $T_h$ satisfies mesh condition \rf{3.10}, then the following superconvergence estimate holds}
\begin{equation}
\|\partial_\beta u-R_h(\partial_\beta u)\|\leq Ch^{k+1}|u|_{k+1},\;k\geq 0.\la{3.20}
\end{equation}
\zm\quad From equations \rf{1.1} and \rf{3.18}, we have
$$
\partial_\beta u-R_h(\partial_\beta u)=-(\alpha +\nabla\cdot\beta)(u-u_h^0).
$$
Hence, estimates \rf{3.19} and \rf{3.20} can be obtained by using \rf{3.9} and \rf{3.16}, respectively, under different mesh conditions.\zb

Compared with the equation to solve numerically the \textcolor[rgb]{.00,0.00,0.00}{directional derivative} $\partial_\beta u$ in the DG method \cite{Cockburn1},   our recovery formula $R_h(\partial_\beta u)$ is very simple and it can be computed explicitly, element by element.

\section{Numerical example}
\setcounter{section}{5} \setcounter{equation}{0}

In this section, we present some numerical examples to show the effectiveness of this WG method.

{\em Example 1}\quad Meshes satisfying condition \rf{3.10}

\begin{figure}[h]
 \begin{center}\begin{picture}(310,100 )(   0.0000000     ,   0.0000000     )

  \def\msq{\begin{picture}(   100.00000     ,   100.00000     )(   0.0000000     ,   0.0000000  )
     \multiput(0,0)(100,0){2}{\line(0,1){100}}\put(100,0){\line(-1,1){100}}
     \multiput(0,0)(0,100){2}{\line(1,0){100}}  \end{picture}  }

\put(0,0){\setlength{\unitlength}{1pt}\msq}
\multiput(110,0)(50,0){2}{\multiput(0,0)(0,50){2}{\setlength{\unitlength}{0.5pt}\msq}}
\multiput(220,0)(25,0){4}{\multiput(0,0)(0,25){4}{\setlength{\unitlength}{0.25pt}\msq}}
 \end{picture}
\end{center}
{\bf Fig. 3}\quad The first three levels of grids for Example 1.
\end{figure}

\begin{figure}[ht] \begin{center}
\begin{picture}(250,100)(0,0)
 \put(0,   0){\scalebox{0.4}{\includegraphics{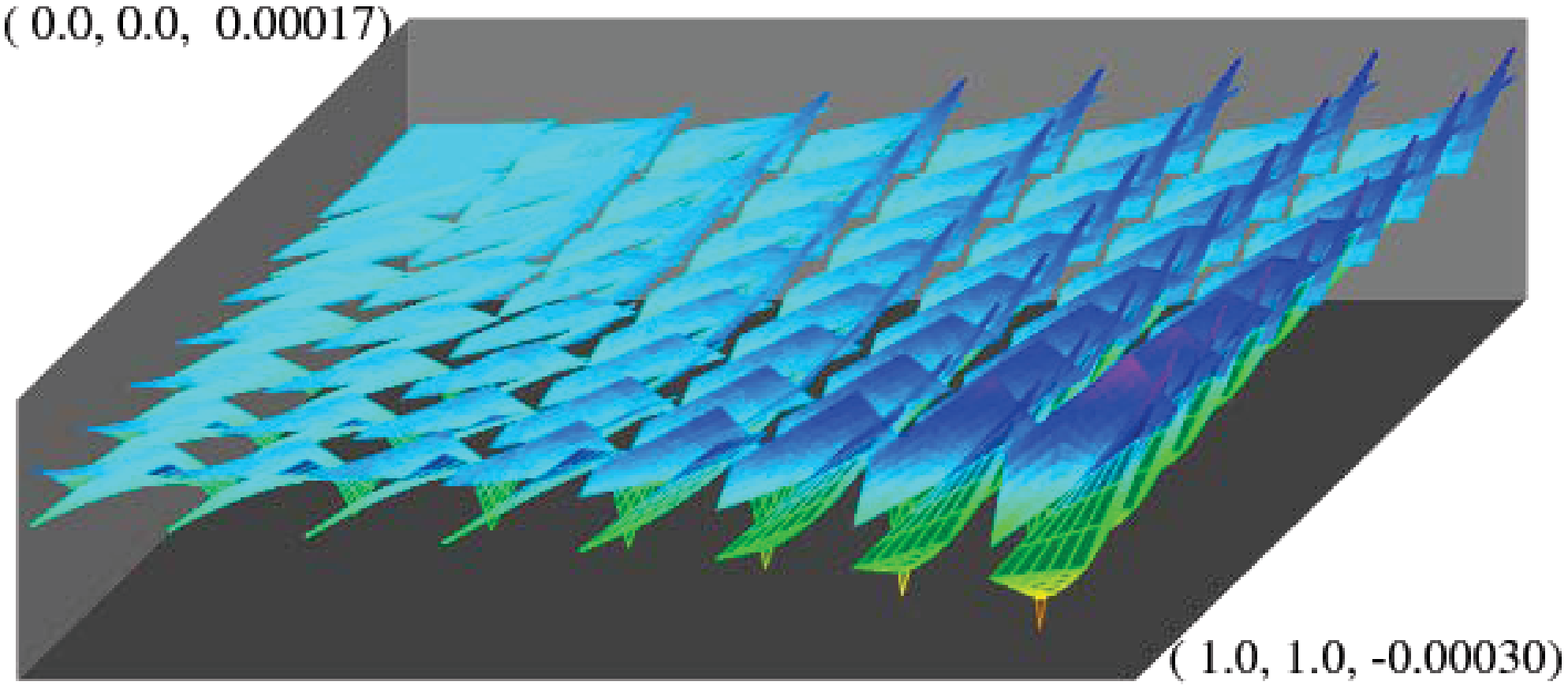}}}
\end{picture}\end{center}
{\bf Fig. 4}\quad The error of $P_2$ WG on the 4th level mesh
    for Example 1.
\end{figure}

\begin{table}[ht]
  \caption{\label{b-n1} The error profile for Example 1
    on meshes of Fig. 3.}
\begin{center}  \begin{tabular}{c|rr|rr|rr}
\hline level  & $ \|u - u^0_h\|$ &rate &
    $ |||Q_h u - u_h|||$ & rate   &
    $ \|\partial_{\beta} u - R_h(\partial_{\beta} u_h)\|$ & rate  \\ \hline
  &\multicolumn{6}{c}{ The $P_1$ WG method  } \\ \hline
 4&   0.1416E-02 & 1.94&   0.6734E-02 & 1.59&   0.2832E-02 & 1.94\\
 5&   0.3618E-03 & 1.97&   0.2300E-02 & 1.55&   0.7236E-03 & 1.97\\
 6&   0.9143E-04 & 1.98&   0.7983E-03 & 1.53&   0.1829E-03 & 1.98\\
  \hline
  &\multicolumn{6}{c}{ The $P_2$ WG method  } \\ \hline
 4&   0.3247E-04 & 2.95&   0.1286E-03 & 2.62&   0.6494E-04 & 2.95\\
 5&   0.4118E-05 & 2.98&   0.2166E-04 & 2.57&   0.8235E-05 & 2.98\\
 6&   0.5181E-06 & 2.99&   0.3729E-05 & 2.54&   0.1036E-05 & 2.99\\
\hline
  &\multicolumn{6}{c}{ The $P_3$ WG method } \\ \hline
 4&   0.5737E-06 & 3.97&   0.2731E-05 & 3.55&   0.1147E-05 & 3.97\\
 5&   0.3617E-07 & 3.99&   0.2366E-06 & 3.53&   0.7233E-07 & 3.99\\
 6&   0.2269E-08 & 3.99&   0.2070E-07 & 3.51&   0.4539E-08 & 3.99\\
\hline
  &\multicolumn{6}{c}{ The $P_4$ WG method } \\ \hline
 3&   0.2601E-06 & 4.76&   0.1161E-05 & 4.59&   0.5203E-06 & 4.76\\
 4&   0.9346E-08 & 4.80&   0.4882E-07 & 4.57&   0.1869E-07 & 4.80\\
 5&   0.3265E-09 & 4.84&   0.2115E-08 & 4.53&   0.6530E-09 & 4.84\\
\hline
\end{tabular}\end{center} \end{table}

We solve problem \rf{1.1} with data (and a exact solution):
\begin{align}\label{son1}
   u=e^{xy}, \quad \beta=(1,0), \quad \alpha=2, \quad \Omega=(0,1)\times(0,1).
\end{align} The computational meshes are displayed in Fig. 3.
Only in this example (not in the next three examples), $\beta \cdot {\bf n}\equiv 0$ on some edges.
\commentout{In this case, we determine the out-flow edge by counter-clock direction of the element boundary,
  and define $\beta\cdot {\bf n}=10^{-13}$ (the machine epsilon) in the computation of
  $s(\cdot,\cdot)$.}
The errors and the orders of convergence are listed in Table 1.
From the table, we can see,  all the WG solutions reach $O(h^{k+1})$-order convergence in
  $L_2$-norm, and $O(h^{k+1/2})$-order convergence in the triple bar norm. For the recovered directional derivative, we get one superconvergence of $O(h^{k+1})$-order in the $L_2$-norm. These numerical results verify the theoretical predictions given in Theorem 3.1, Theorem 3.2 and Theorem 4.1.

To understand error behavior,  we plot the error of $P_2$ WG solution on the 4th level mesh in
  Fig. 4.

{\em Example 2}\quad Non-compatible meshes

In this example, we solve problem \rf{1.1} with data (and a solution):
\begin{align}\label{so2}
   u=\sin 4x \sin 4y, \quad \beta=(1,1), \quad \alpha=1, \quad \Omega=(0,1)\times(0,1).
\end{align} The computation meshes are displayed in Fig.5.
The errors and the orders of convergence are listed in Table \ref{b-2}.
From this table, we can see that all the WG solutions reach $O(h^{k+1})$-order convergence in
the $L_2$-norm, and $O(h^{k+1/2})$-order convergence in the triple bar norm.
For the recovered directional derivative,  we get one superconvergence of
   $O(h^{k+1})$-order in the $L_2$-norm. Note that in this example, the mesh is non-compatible and condition \rf{3.10} is not satisfied, but we still obtain the $O(h^{k+1})$-order convergence. What is the necessary condition for the meshes to produce the (k+1)-order convergence, this is an open question.
\begin{figure}[h]
 \begin{center}\begin{picture}(170,80 )(   0.0000000     ,   0.0000000     )
\put(0,0){\setlength{\unitlength}{0.4pt}\meig}
\multiput(85,0)(40,0){2}{\multiput(0,0)(0,40){2}{\setlength{\unitlength}{0.2pt}\meig}}
 \end{picture}\end{center}
\centerline{{\bf Fig. 5} The first two levels of meshes for Example 2. }
\end{figure}

\begin{table}[ht]
  \caption{\label{b-2} The error profile for Example 2 on meshes of Fig. 6.}
\begin{center}  \begin{tabular}{c|rr|rr|rr}
\hline level  & $ \|u - u^0_h\|$ &rate &
    $ |||Q_h u - u_h|||$ & rate   &
    $ \|\partial_{\beta} u - R_h(\partial_{\beta} u_h)\|$ & rate  \\ \hline
  &\multicolumn{6}{c}{ The $P_1$ WG method } \\ \hline
 5&   0.2149E-02 & 2.00&   0.2042E-01 & 1.50&   0.2149E-02 & 2.00\\
 6&   0.5372E-03 & 2.00&   0.7210E-02 & 1.50&   0.5372E-03 & 2.00\\
 7&   0.1344E-03 & 2.00&   0.2547E-02 & 1.50&   0.1344E-03 & 2.00\\
  \hline
  &\multicolumn{6}{c}{ The $P_2$ WG method } \\ \hline
 4&   0.4157E-03 & 3.04&   0.3586E-02 & 2.50&   0.4157E-03 & 3.04\\
 5&   0.5210E-04 & 3.00&   0.6333E-03 & 2.50&   0.5210E-04 & 3.00\\
 6&   0.6574E-05 & 2.99&   0.1119E-03 & 2.50&   0.6574E-05 & 2.99\\
\hline
  &\multicolumn{6}{c}{ The $P_3$ WG method } \\ \hline
 4&   0.2472E-04 & 4.04&   0.2438E-03 & 3.51&   0.2472E-04 & 4.04\\
 5&   0.1565E-05 & 3.98&   0.2149E-04 & 3.50&   0.1565E-05 & 3.98\\
 6&   0.1015E-06 & 3.95&   0.1896E-05 & 3.50&   0.1015E-06 & 3.95\\
\hline
  &\multicolumn{6}{c}{ The $P_4$ WG method } \\ \hline
 3&   0.4561E-04 & 5.29&   0.3570E-03 & 4.54&   0.4561E-04 & 5.29\\
 4&   0.1373E-05 & 5.05&   0.1557E-04 & 4.52&   0.1373E-05 & 5.05\\
 5&   0.4249E-07 & 5.01&   0.6834E-06 & 4.51&   0.4249E-07 & 5.01\\
\hline
\end{tabular}\end{center} \end{table}

\begin{figure}[ht]
 \begin{center}\begin{picture}( 170,80 )(   0.0000000     ,   0.0000000     )

\put(0,0){\setlength{\unitlength}{0.4pt}\msev}
\multiput(85,0)(40,0){2}{\multiput(0,0)(0,40){2}{\setlength{\unitlength}{0.2pt}\msev}}
 \end{picture}\end{center}
\centerline{{\bf Fig. 6}. The first two levels of grids for Example 3}
\end{figure}

 \begin{table}[ht]
  \caption{\label{b-3} The error profile for Example 3 on meshes of Fig. 7.}
\begin{center}  \begin{tabular}{c|rr|rr|rr}
\hline level  & $ \|u - u^0_h\|$ &rate&
    $ |||Q_h u - u_h|||$ & rate   &
    $ \|\partial_{\beta} u - R_h(\partial_{\beta} u_h)\|$ & rate  \\ \hline
  &\multicolumn{6}{c}{ The $P_1$ WG method} \\ \hline
 5&   0.1697E-02 & 1.94&   0.2357E-01 & 1.47&   0.5092E-02 & 1.94\\
 6&   0.4357E-03 & 1.96&   0.8411E-02 & 1.49&   0.1307E-02 & 1.96\\
 7&   0.1117E-03 & 1.96&   0.2987E-02 & 1.49&   0.3350E-03 & 1.96\\
  \hline
  &\multicolumn{6}{c}{ The $P_2$ WG method } \\ \hline
 3&   0.1864E-02 & 3.05&   0.1647E-01 & 2.53&   0.5592E-02 & 3.05\\
 4&   0.2396E-03 & 2.96&   0.2864E-02 & 2.52&   0.7189E-03 & 2.96\\
 5&   0.3288E-04 & 2.87&   0.5010E-03 & 2.52&   0.9863E-04 & 2.87\\\hline
  &\multicolumn{6}{c}{ The $P_3$ WG method } \\ \hline
 3&   0.9985E-04 & 4.01&   0.7953E-03 & 3.49&   0.2995E-03 & 4.01\\
 4&   0.6304E-05 & 3.99&   0.7069E-04 & 3.49&   0.1891E-04 & 3.99\\
 5&   0.4179E-06 & 3.91&   0.6267E-05 & 3.50&   0.1254E-05 & 3.91\\\hline
\end{tabular}\end{center} \end{table}

{\em Example 3}\quad \textcolor[rgb]{.00,0.00,0.00}{Non-divergence-free flow and non-compatible meshes}

In this example, we solve problem \rf{1.1} with data (and a solution):
\begin{align}\label{so3}
   u=(x+y)^2(x+y-1)^2, \quad \beta=(x,y), \quad \alpha=1, \quad \Omega=(0,1)\times(0,1).
\end{align} The computation meshes are displayed in Fig. 6.
The errors and the orders of convergence are listed in Table \ref{b-3}.
From this table, we can see that all the WG solutions reach $O(h^{k+1})$-order convergence in
 the $L_2$-norm, and $O(h^{k+1/2})$-order convergence in the triple bar norm.
For the recovered directional derivative, we get one superconvergence of
   $O(h^{k+1})$-order in the $L_2$-norm.

{\em Example 4}\quad A circular flow

We solve problem \rf{1.1} with data:
\begin{align}\label{so4}
   \beta&=(-y,x), \quad \alpha=0, \quad f=0, \quad \Omega=(-1,1)^2\setminus [0,1]\times\{0\}, \\
       g&=\begin{cases} \sin^2\pi x, & \hbox{ on the inflow boundary } [0,1]\times\{0^+\}, \\
                      0, & \hbox{ on the rest inflow boundary } \{0^-\}\times [0,1], [-1,0]\times\{1^-\},\\
                  &\qquad\qquad \{-1^-\}\times [-1,0], [0,1]\times\{1^+\}. \end{cases} \nonumber
\end{align}
The computation meshes are displayed in Fig. 7.
The exact solution of \eqref{so4} is unknown.
But as the inflow is minimum,  we expect that the final outflow profile is very close
  to that of the inflow.
This is observed in the numerical solutions, plotted in Fig. 8.
\begin{figure}[h]
 \begin{center}\begin{picture}(310,100 )(   0.0000000     ,   0.0000000     )

  \def\msq{\begin{picture}(   100.00000     ,   100.00000     )(   0.0000000     ,   0.0000000  )
     \multiput(0,0)(100,0){2}{\line(0,1){100}}
     \multiput(0,0)(0,100){2}{\line(1,0){100}}  \end{picture}  }

\put(0,0){\setlength{\unitlength}{1pt}\msq}
\multiput(110,0)(50,0){2}{\multiput(0,0)(0,50){2}{\setlength{\unitlength}{0.5pt}\msq}}
\multiput(220,0)(25,0){4}{\multiput(0,0)(0,25){4}{\setlength{\unitlength}{0.25pt}\msq}}
 \end{picture}\end{center}
\centerline{{\bf Fig. 7.} The first three levels of grids for Example 4.}
\end{figure}

\begin{center}
 \scalebox{0.4}{\includegraphics{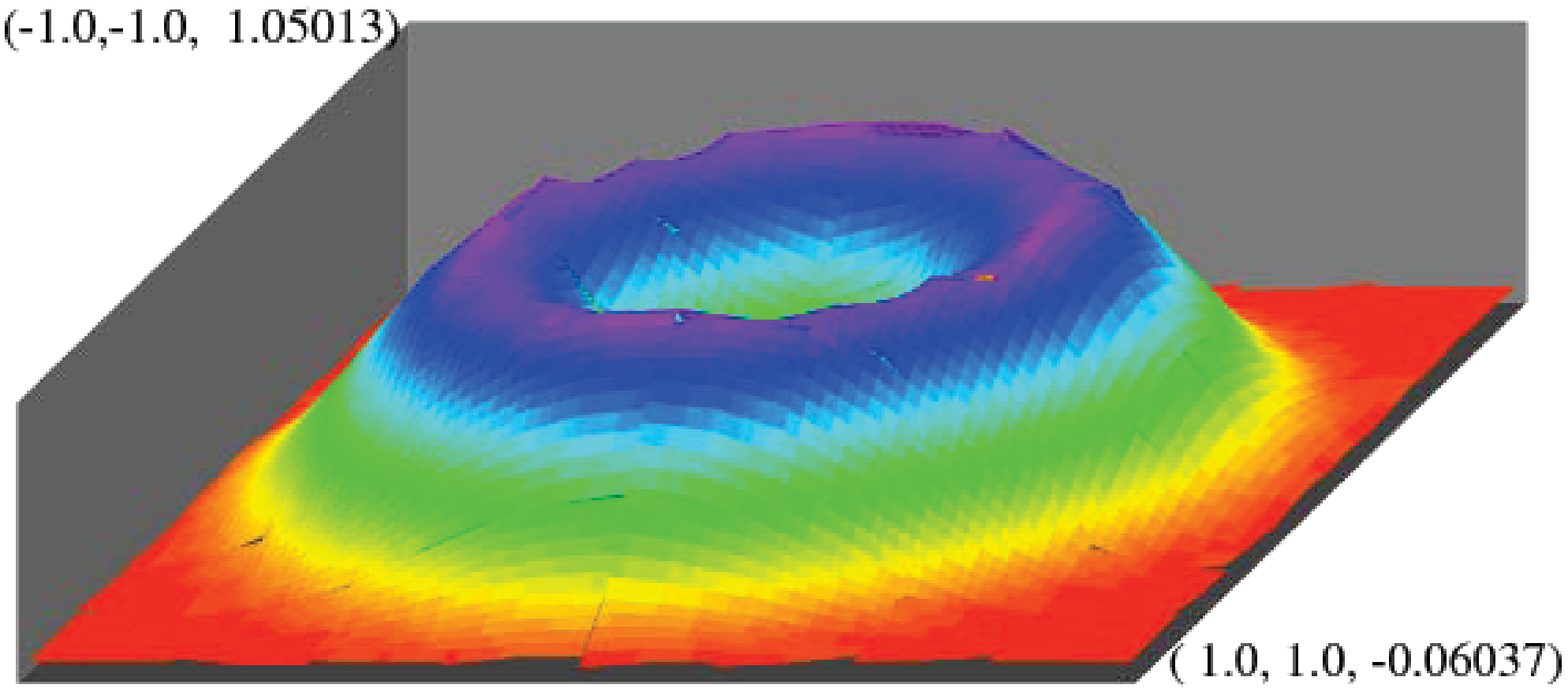}}
 \scalebox{0.4}{\includegraphics{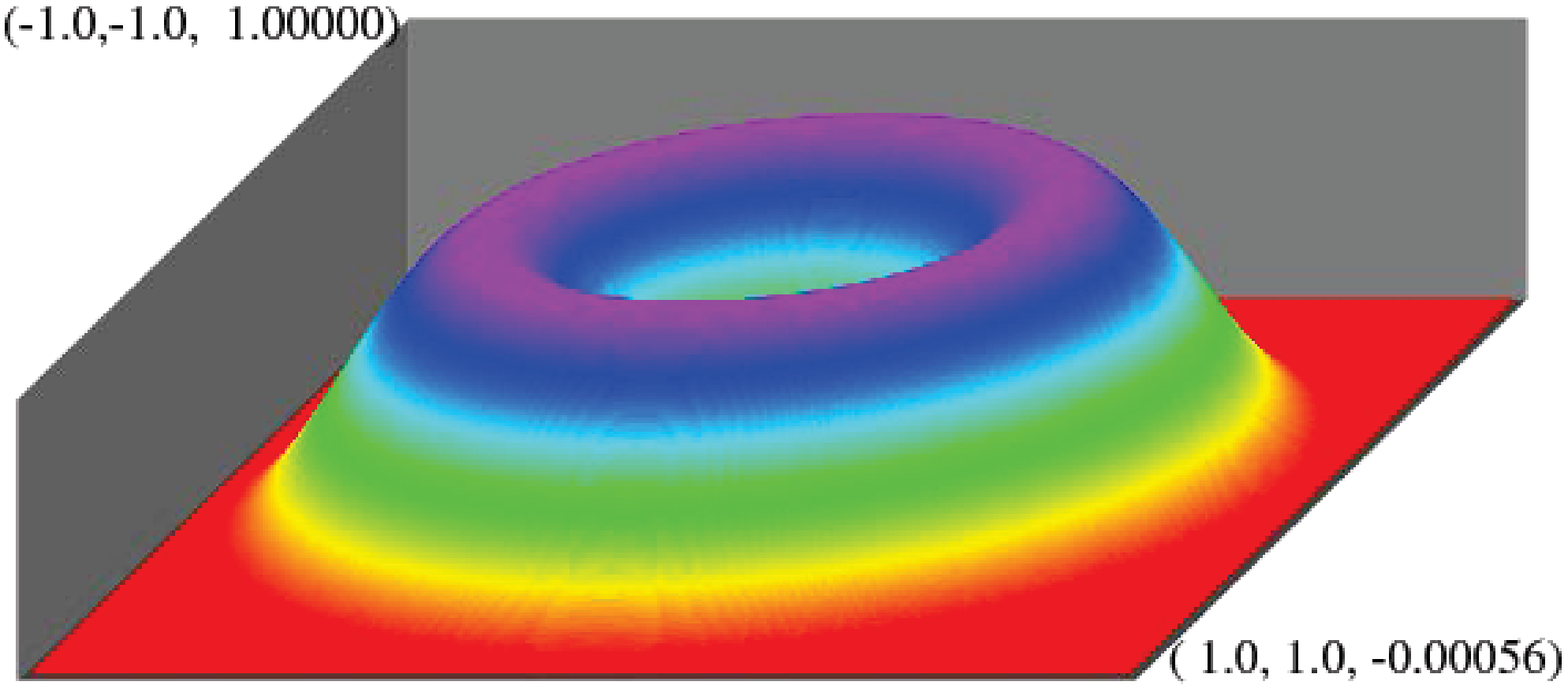}}
\end{center}

{\bf Fig. 8.} The numerical solution of $P_2$ WG on the 3rd level mesh, and the $P_4$ WG solution
    on the 5th level mesh for Example 4.
\section{Conclusion}
We present and analyze a weak Galerkin finite element method for solving the typical first order hyperbolic equation: the transport-reaction equation. This method is highly flexible by allowing the use of discontinuous finite element on general meshes consisting of arbitrary polygon/polyhedra. \commentout{Like the original DG method for this problem, the discrete WG equation can be solved locally, element by element.} Using the $k$th-order polynomials ($k\geq 0$), we prove that the WG solution admits the optimal $L_2$-convergence rate of $O(h^{k+1})$-order under special mesh condition which is slightly weaker than the {\em flow condition} given in \cite{Cockburn1} for the DG method solving this problem.
Moreover, a derivative recovery formula is presented to approximate the convection \textcolor[rgb]{.00,0.00,0.00}{directional derivative} and the corresponding superconvergence estimate is given. Numerical examples on compatible and non-compatible meshes are provided to show the effectiveness of this WG method. Our work provides an approach to develop the WG method for first order hyperbolic problems.

\section*{Acknowledgments}
This work was supported by the State Key Laboratory of Synthetical
Automation for Process Industries Fundamental Research Funds, No. 2013ZCX02.

\end{document}